\documentclass{article}
\usepackage{graphicx,amssymb,mathrsfs,amsmath}
\newtheorem{theorem}{Theorem}

\newtheorem{nota}{Notation}
\newtheorem{definition}{Definition}
\newcommand{\qed}{\hfill\rule{0.3ex}{1.5ex}}
\begin{document}
 \title{An equation for the general Ramsey number $R(P_1,P_2,...,P_t;r)$} \author{Kunjun Song\\
Department of Modern Physics\\ University of Science and Technology of China\\ Hefei,Anhui,P.R.China\\
\texttt{skjmom@mail.ustc.edu.cn}}
 \date{}
 \maketitle
 \begin{abstract}
  The Ramsey number $R(P_1,P_2,...,P_t;r)$ is a valve value such that as
long as the cardinality $n$ of the $n$-set $V_n=\{1,\cdots,n\}$ is no less than $R$,however all the $\binom{n}{r}$ $r$-subsets of $V_n$ are distributed
into $t$ boxes, $V_n$ will always have a property $W$ expressed as eq.(1).Thus, by calculating the number of ways of distribution of
$r$-subsets that makes $W$ true,one can get an equation for $R(P_1,P_2,...,P_t;r)$.The evaluation of the general term in this eq.
and the counting of the frequencies of occurrence of the various values the general term takes can be reduced to the problem of
elementary counting.
 \end{abstract}

Roughly speaking,Ramsey theory is the precise mathematical formulation of the statement:\emph{Complete disorder is impossible.} or \emph{Every large enough structure will inevitably contain some regular substructures}. The Ramsey number measures how large on earth does the structure need to be so that the specified substructures are guaranteed to emerge.

The most general (finite) Ramsey number[1] is defined by the following existence theorem:

 $R(P_1,P_2,...,P_t;r)$ is the smallest integer $n$ that has the following property.All the $t^{\binom{n}{r}}$ ways of distribution
 of the $\binom{n}{r}$ $r$-subsets of the $n$-set $V_n=\{1,\cdots,n\}$ into $t$ boxes makes the following event $W$  true: There exists a $P_1$-subset,all $r$-subsets of which are in box 1;or there exists a $P_2$
-subset,all $r$-subsets of which are in box 2;or...or there exists a $P_t$-subset,all $r$-subsets of which are in box $t$.

 This definition is equivalent to saying that the Ramsey number \\$R(P_1,P_2,...,P_t;r)$ is the smallest positive integer
 $n$ that satisfies the following equation:(where we have use the principle of inclusion and exclusion)\footnote{Here I assume some method has been devised to order the various
$P_i$-events.For example,$A_{i_1j_1}<A_{i_2j_2}$ iff $i_1<i_2$ or $i_1=i_2$ and $j_1<j_2$.c.f.appendix 1.}
\begin{align}
 &N(W)\equiv N(\bigcup_{i=1}^{t}\bigcup_{j=1}^{\binom{n}{P_i}}A_{ij})\notag\\
&=\sum_{k=1}^{\binom{n}{P_1}+\cdots+\binom{n}{P_t}}(-1)^{k-1}\sum_{(1,1)\leq(i_1,j_1)<\cdots<(i_k,j_k)\leq(t,\binom{n}{P_t})}N(A_{i_1j_1}\bigcap\cdots\bigcap
A_{i_kj_k})\notag\\ &=t^{\binom{n}{r}}\equiv \mbox{txp}\left[\binom{n}{r}\right]
\end{align}

Here we have introduced the following notations:
\begin{nota}
 $A_{ij}$ denotes the event ``All the $\binom{P_i}{r}$ $r$-subsets of
the $j^{th}$ $P_i$ -subset of the $n$-set $V_n$ are in the $i^{th}$ box\footnote{See appendix 1 for the order of the $\binom{n}{P_i}$ $P_i$-subsets.}." The corresponding $j^{th}$ $P_i$ -subset is denoted by
$\bar{A}_{ij}$.($1\leq i\leq t,1\leq j\leq \binom{n}{P_i};r\le P_1\le P_2\le\cdots\le P_t$) \qed
\end{nota}
\begin{nota}
The number
of ways of distribution of $r$-subsets into $t$  boxes that makes event $X$ true is denoted by $N(X)$.\qed
\end{nota}

We remark that it is the event $A_{ij}$,not the subset $\bar{A}_{ij}$, of importance,as it is possible when $i_\mu\neq
i_\nu$,$P_{i_\mu}=P_{i_\nu}$.When this happens, $\bar{A}_{i_{\mu}j}=\bar{A}_{i_{\nu}j}$,but $A_{i_{\mu}j}\neq A_{i_{\nu}j}$ .  For
future purpose,we need a few more definitions:
\begin{definition}
The collection $\bigcup_{i=1}^{t}\bigcup_{j=1}^{\binom{n}{P_i}}A_{ij}$ of all the $P_1$-events,$P_2$-events,...,\-$P_t$-events is denoted by $V_{P_1,\cdots,P_t}$.The element $(A_{i_1j_1},\cdots,A_{i_kj_k})$;$A_{i_1j_1}<\cdots<$\\$A_{i_kj_k}$ in the Cartesian product
$V_{P_1,\cdots,P_t}^k$ \-is called an unordered $k$-event tuple.The corresponding element $(\bar{A}_{i_1j_1},\cdots,\bar{A}_{i_kj_k})\in \bar{V}_{P_1,\cdots,P_t}^k\equiv\bigcup_{i=1}^{t}\bigcup_{j=1}^{\binom{n}{P_i}}\bar{A}_{ij}$ is called an unordered $k$-set tuple.\qed
\end{definition}
\begin{nota}
The map that maps an unordered $s$-set tuple to the number of elements in the intersection of the $s$ sets\footnote{These sets are not necessarily distinct,i.e it is possible that there are actually less than $s$ sets,but there will
always be $s$ distinct corresponding events,and the stated degenerate case is of no interest as it is not compatible in the sense of
Def.2.} $\bar{A}_{i_1j_1}\cdots \bar{A}_{i_sj_s}$ is denoted by\footnote{When the argument(the unordered $s$-set tuple $(\bar{A}_{i_1j_1},\cdots, \bar{A}_{i_sj_s})$ ) of the function $P_{i_1\cdots i_s}$ is clear from the context,it will usually be omitted.}
\begin{align}
&P_{i_1\cdots i_s}:\bar{V}_{P_{i_1},\cdots,P_{i_s}}^s\rightarrow \mathbb{N}\notag\\
&P_{i_1\cdots i_s}(\bar{A}_{i_1j_1},\cdots, \bar{A}_{i_sj_s})\equiv
P_{i_1\cdots i_s}(j_1\cdots j_s)=\left|\bigcap_{m=1}^{s}A_{i_mj_m}\right|\notag
\end{align}
\qed
\end{nota}

\begin{definition}
 The unordered $k$-event tuple $(A_{i_1j_1},\cdots,A_{i_kj_k})$ is called distributionally compatible iff for any two events $A_{i_\mu
j_\mu}$ and $A_{i_\nu j_\nu}$ chosen arbitrarily from the $k$ events $A_{i_1j_1}\cdots A_{i_kj_k}$, the relation $P_{i_\mu
i_\nu}\le r-1$ holds.($\forall 1\le\mu<\nu\le k$ ,i.e. $\bar{A}_{i_\mu j_\mu}$ and $\bar{A}_{i_\nu j_\nu}$ must actually have no common
$r$-subsets if the two events corresponding to them demand their common $r$-subsets be put into two different boxes $i_\mu\neq i_\nu$.)\footnote{When
$i_\mu=i_\nu$,i.e. if  the two events correspond to them demand their common $r$-subsets be put into the same box,then naturally we have
$P_{i_\mu i_\nu}\le P_{i_\mu}$.This is not a constraint on the value of $P_{i_\mu i_\nu}$.}\qed
\end{definition}
\begin{definition}
In the Venn diagram of the $k$ sets $\bar{A}_{i_1j_1}\cdots\bar{A}_{i_kj_k}$,$\bar{A}_{i_1j_1}<\cdots<\bar{A}_{i_kj_k}$ the
$n$-set $V_n$ is divided into $2^k$ disjoint parts,which will be termed the Venn parts of $V_n$ w.r.t. these $k$ sets.

Let's use the $2^k$ $k$-digit binary numbers to represent the $2^k$ Venn parts.The $m^{th}$ digit of any one of these numbers being
0 or 1 depends on whether the corresponding part is contained(1) in $A_{i_mj_m}$ or not(0).The cardinality of the Venn part represented
by the binary number $B$ will be denoted by $Q^{(k)}_B(\bar{A}_{i_1j_1},\cdots,\bar{A}_{i_kj_k})$.\footnote{Similar to footnote 4,when the argument $(\bar{A}_{i_1j_1},\cdots,\bar{A}_{i_kj_k})$ of the function $Q^{(k)}_B$ is understood,its value $Q^{(k)}_B(\bar{A}_{i_1j_1},\cdots,\bar{A}_{i_kj_k})$ will simply be denoted by $Q^{(k)}_B$.} \qed
\end{definition}

\begin{definition}
The map
 \begin{align}
V_{P_1,\cdots,P_t}^k\ni(A_{i_1j_1},\cdots,A_{i_kj_k})\mapsto(Q^{(k)}_{B_1},Q^{(k)}_{B_2},\cdots,Q^{(k)}_{B_{2^k}})\in
\mathbb{N}^{2^k}
 \end{align}
  where $\bar{A}_{i_1j_1}<\cdots<\bar{A}_{i_kj_k}$ and $B_1<B_2<\cdots<B_{2^k}$,is called the Venn spectrum of the unordered $k$-event tuple
$(A_{i_1j_1},\cdots,A_{i_kj_k})$.\qed
\end{definition}

\begin{definition}
The map
\begin{align}
&V_{P_1,\cdots,P_t}^k\ni(A_{i_1j_1},\cdots,A_{i_kj_k})\mapsto\\
&(n,\{P_{i_\lambda}\}_{\lambda=1}^k,\{P_{i_{\lambda_1}i_{\lambda_2}}\}_{1\le\lambda_1<\lambda_2\le
k},\cdots,\notag\{P_{i_{\lambda_1}\cdots i_{\lambda_s}}\}_{1\le\lambda_1<\cdots<\lambda_s\le k},\cdots,P_{i_1\cdots i_k})\in\mathbb{N}^{2^k}
 \end{align}
 where $\bar{A}_{i_1j_1}<\cdots<\bar{A}_{i_kj_k}$ and $i_{\lambda_1}<\cdots<i_{\lambda_s}$,is called the intersection spectrum of the
unordered $k$-event tuple $(A_{i_1j_1},\cdots,A_{i_kj_k})$.\qed
 \end{definition}

We remark that it is the $Q$'s that will be used as basic variables in the following.Any $P$'s that appear below should be
understood to be an abbreviation of a sum of the $Q$'s according to Theorem 1.\footnote{According to Theorem 1,to get the $P$'s from the $Q$'s,we only need to do addition,while the reverse requires inclusion-exclusion.This is one of the reasons why we use the $Q$'s as basic variables.}

\begin{theorem}
The Venn spectrum and intersection spectrum of an unordered $k$-event tuple $(A_{i_1j_1},\cdots,A_{i_kj_k})$ is connected by the following linear equations:
\begin{align}
 P_{i_{\lambda_1}\cdots i_{\lambda_s}}=\sum\mbox{all $2^{k-s}$ Q's whose $\lambda_1^{th}\cdots\lambda_s^{th}$
digits are 1}
 \end{align}
 \begin{equation} Q^{(k)}_{11\cdots 1}=P_{i_1i_2\cdots i_k}\notag
 \end{equation}
  \begin{equation}
Q^{(k)}_{01\cdots 1}=P_{i_2\cdots i_k}-P_{i_1\cdots i_k},\cdots, Q^{(k)}_{1\cdots 10}=P_{i_1i_2\cdots
i_{k-1}}-P_{i_1\cdots i_k}\notag
 \end{equation}
  \begin{align}
  Q^{(k)}_{(\mbox{$k$-digit binary number:
$\nu_1^{th}\cdots\nu_{k-f}^{th}$ digits 1,other $f$ digits 0)}}\notag
\end{align}
 \begin{equation}
=P_{i_{\nu_1}i_{\nu_2}\cdots i_{\nu_{,k-f}}}+\sum_{s=1}^f(-)^s\sum_{\begin{subarray}{c}
1\leq\lambda_1<\lambda_2<\cdots<\lambda_s\leq k\\
\lambda_m\neq\nu_q (\forall m=1\cdots s,\forall q=1\cdots k-f)
\end{subarray}
}P_{i_{\nu_1}i_{\nu_2}\cdots
i_{\nu_{,k-f}}i_{\lambda_1}i_{\lambda_2}\cdots i_{\lambda_s}}\notag
 \end{equation}
  \begin{equation}
   Q^{(k)}_{00\cdots
0}=n+\sum_{s=1}^k(-)^s\sum_{1\leq\lambda_1<\lambda_2<\cdots<\lambda_s\leq k}P_{i_{\lambda_1}i_{\lambda_2}\cdots i_{\lambda_s}}
\end{equation}\qed
\end{theorem}
 \textbf{Proof}:
  Given the $Q$'s,it is really a matter of inspection and induction with the help of
Venn diagram to get the $P$'s, which (eq.(4)) can also be verified immediately by a little thought of the very meaning of the $Q$'s
and the $P$'s,It requires more work (inclusion-exclusion) to get the $Q$'s from the $P$'s,but this can also be done by inspection of some simple cases and then generalize,as this result will not be used often in the following,we leave its proof as an exercise.
$\hfill\blacksquare$

In order to bring the basic eq.(1) into a more explicit and convenient form,we first write it schematically as
\begin{align}
&\sum_{k=1}^{\binom{n}{P_1}+\cdots+\binom{n}{P_t}}(-)^{k-1}\sum_{\mbox{value}}\mbox{value of $N(A_{i_1j_1}\bigcap\cdots\bigcap
A_{i_kj_k})$}\times\mbox{frequency this value occurs}\notag\\ &=t^{\binom{n}{r}}
 \end{align}

We can now make our claim about the first factor in the second summation $\sum_{\mbox{value}}$ in eq.(6).\\
\begin{theorem}
 If the unordered $k$-event tuple $(A_{i_1j_1},\cdots,A_{i_kj_k})$ is distributionally compatible,then
  \begin{align}
   N(A_{i_1j_1}\bigcap\cdots\bigcap A_{i_kj_k})
&=\mbox{txp}\left[\binom{n}{r}+\sum_{s=1}^{k}(-1)^s\sum_{1\leq\lambda_1<\cdots<\lambda_s\leq k}\binom{P_{i_{\lambda_1}\cdots
i_{\lambda_s}}}{r}\right]
\end{align}
,otherwise
\begin{align}
 N(A_{i_1j_1}\bigcap\cdots\bigcap A_{i_kj_k})=0
 \end{align}\qed
\end{theorem}
\textbf{Proof}:
 The event $A_{i_1j_1}\bigcap\cdots\bigcap A_{i_kj_k}$ demands the
$\left|\bigcup_{s=1}^{k}\mbox{Tr}(\bar{A}_{i_sj_s})\right|$ $r$-subsets (Here $\mbox{Tr}(X)$ represents the family consisting of all the
$r$-subsets of the set $X$) be simultaneously put into the appropriate boxes specified by this event.When the unordered $k$-event tuple
$(A_{i_1j_1},\cdots,A_{i_kj_k})$ is distributionally compatible,this can be done and after doing this,the remaining
$\binom{n}{r}-\left|\bigcup_{s=1}^{k}\mbox{Tr}(\bar{A}_{i_sj_s})\right|$ $r$-subsets can then be arbitrarily put into the $t$ boxes,thus
$N(A_{i_1j_1}\bigcap\cdots\bigcap
A_{i_kj_k})=$\\$\mbox{txp}\left[\binom{n}{r}-
\left|\bigcup_{s=1}^{k}\mbox{Tr}(\bar{A}_{i_sj_s})\right|\right]$ .Then we expand
$\left|\bigcup_{s=1}^{k}\mbox{Tr}(\bar{A}_{i_sj_s})\right|$ using inclusion-exclusion,
\begin{align}
 N(A_{i_1j_1}\bigcap\cdots\bigcap
A_{i_kj_k})&=\mbox{\mbox{txp}}\left[\binom{n}{r}-\left|\bigcup_{\alpha=1}^{k}\mbox{Tr}(\bar{A}_{i_\alpha j_\alpha})\right|\right]\\
&=\mbox{\mbox{txp}}\left[\binom{n}{r}+\sum_{s=1}^{k}(-1)^s\mathop{\sum}_{1\leq\lambda_1<\cdots<\lambda_s\leq
k}\left|\bigcap_{m=1}^{s}\mbox{Tr}(\bar{A}_{i_{\lambda_m}j_{\lambda_m}})\right|\right]\notag\\
&=\mbox{\mbox{txp}}\left[\binom{n}{r}+\sum_{s=1}^{k}(-1)^s\sum_{1\leq\lambda_1<\cdots<\lambda_s\leq k}\binom{P_{i_{\lambda_1}\cdots
i_{\lambda_s}}}{r}\right]\notag
\end{align}
 where use have been made of the fact
  \begin{align}
\left|\bigcap_{m=1}^{s}\mbox{Tr}(\bar{A}_{i_{\lambda_m}j_{\lambda_m}})\right|=
\left|\mbox{Tr}(\bigcap_{m=1}^{s}\bar{A}_{i_{\lambda_m}j_{\lambda_m}})\right|=
\binom{\left|\bigcap_{m=1}^{s}\bar{A}_{i_{\lambda_m}j_{\lambda_m}}\right|}{r} \equiv\binom{P_{i_{\lambda_1}i_{\lambda_2}\cdots
i_{\lambda_s}}}{r}\notag
\end{align}
This in turn comes from
$\bigcap_{m=1}^{s}\mbox{Tr}(\bar{A}_{i_{\lambda_m}j_{\lambda_m}})=\mbox{Tr}(\bigcap_{m=1}^{s}\bar{A}_{i_{\lambda_m}j_{\lambda_m}})$. i.e. the
common $r$-subsets of the $s$ sets $\bar{A}_{i_{\lambda_m}j_{\lambda_m}} \quad(1\leq m\leq s)$ can only be the $r$-subsets of their
intersection,since the elements of the common $r$-subsets must be contained in all the $s$ sets,therefore belonging to their
intersection.The second part of the Theorem is trivial,since when the unordered $k$-event tuple $(A_{i_1j_1},\cdots,A_{i_kj_k})$ is not distributionally compatible,these $k$
events cannot be simultaneously true.This is because when $i_\mu\neq i_\nu$ and $P_{i_\mu i_\nu}> r-1$, the two sets $\bar{A}_{i_\mu j_\mu}$ and $\bar{A}_{i_\nu j_\nu}$ have common $r$-subsets,and the corresponding events $A_{i_\mu j_\mu}$ and $A_{i_\nu j_\nu}$
demand the common $r$-subsets be put into two different boxes $i_\mu\neq i_\nu$ ,which is impossible and hence
$N(A_{i_1j_1}\bigcap\cdots\bigcap A_{i_kj_k})=0$ in this case.$\hfill\blacksquare$

We now formulate our result about the second factor in the second summation $\sum_{\mbox{value}}$ in eq.(6).\\
\begin{theorem}
The frequency of occurrence of the value \\
$\mbox{\mbox{txp}}\left[\binom{n}{r}+\sum_{s=1}^{k}(-1)^s\sum_{1\leq\lambda_1<\cdots<\lambda_s\leq k}\binom{P_{i_{\lambda_1}\cdots
i_{\lambda_s}}}{r}\right]$ of $N(A_{i_1j_1}\bigcap\cdots\bigcap A_{i_kj_k})$ is given by the Venn spectrum of the unordered $k$-event tuple
$(A_{i_1j_1},\cdots,A_{i_kj_k})$ as follows
\begin{align}
n!/\prod_{\mbox{all $2^k$ $k$-digit binary
numbers $B$ }}Q^{(k)}_B!
\end{align}
\qed
 \end{theorem}
  \textbf{Proof}:
  We want to count how many compatible unordered $k$ tuples are there that correspond to a particular value $\mbox{txp}\left[\binom{n}{r}+\sum_{s=1}^{k}(-1)^s\sum_{1\leq\lambda_1<\cdots<\lambda_s\leq k}\binom{P_{i_{\lambda_1}\cdots
i_{\lambda_s}}}{r}\right]$ of $N(A_{i_1j_1}\bigcap\cdots\bigcap
A_{i_kj_k})$.

To find this number,let's consider the following series of maps:(where \{...\} means the set whose elements are...)
\begin{align}
&\mbox{\{Value\}$\mapsto$\{intersection spectrum\}$\mapsto$\{Venn spectrum\}}\\
&\mbox{$\mapsto$\{ways of distribution of $n$ elements into $2^k$ Venn part}\notag\\
&\mbox{whose cardinalities are specified by the Venn spectrum\}$\mapsto$}\notag\\
&\mbox{\{unordered $k$-set tuple\}$\mapsto$\{distributionally compatible unordered $k$-event tuple\}}\notag
\end{align}

We claim that \\
\begin{itemize}
\item 1 Value corresponds to $k!$ intersection spectrum ;\\
\item 1 intersection spectrum corresponds to 1 Venn spectrum;\\
\item 1  Venn spectrum corresponds to  $n!/\prod_BQ^{(k)}_B!$ ways of distribution;\\
\item $k!$ ways of distribution corresponds to 1 unordered $k$-set tuple; \\
\item 1 unordered $k$-set tuple corresponds to 1 compatible unordered  $k$-event tuple.\\
\end{itemize}

This claim obviously implies the Theorem.

Let's now prove this claim.

\begin{itemize}
\item 1 Value corresponds to $k!$ intersection spectrum:\\This is because the value $\mbox{txp}\left[\binom{n}{r}+\sum_{s=1}^{k}(-1)^s\sum_{1\leq\lambda_1<\cdots<\lambda_s\leq k}\binom{P_{i_{\lambda_1}\cdots
i_{\lambda_s}}}{r}\right]$ of $N(A_{i_1j_1}\bigcap\cdots\bigcap
A_{i_kj_k})$ is a funtion of the intersection spectrum $P_{i_{\lambda_1}\cdots
i_{\lambda_s}}$ that is invariant under a permutation of the indices $(i_{\lambda_1}j_{\lambda_1},\cdots,i_{\lambda_s}j_{\lambda_s})$ of $P_{i_{\lambda_1}\cdots
i_{\lambda_s}}(j_{\lambda_1},\cdots,
j_{\lambda_s})$,namely the reordering of the $k$ sets $A_{i_1j_1},\cdots,
A_{i_kj_k}$,which changes the intersection spectrum but does not change the value.

\item 1 intersection spectrum corresponds to 1 Venn spectrum:\\
This is the content of Theorem 1.

\item 1  Venn spectrum corresponds to  $n!/\prod_BQ^{(k)}_B!$ ways of distribution:\\
Note that the $2^k$ Venn parts are ordered (or discernable),then this follows from a classic result in enumerative combinatorics.

\item $k!$ ways of distribution corresponds to 1 unordered $k$-set tuples:\\
For a given way of distribution of $n$ elements into $2^k$ Venn parts that generates a particular
unordered $k$-set tuple.We can perform $k!$ permutations of the $k$ sets,each of which generates a new way of distribution and does not change the particular unordered $k$-set tuple.Each of the $k!$ permutations of the $k$ sets is in turn composed of a number of simultaneous permutations of the corresponding Venn parts contained in the sets being permuted.The net effect of these simultaneous permutations of the corresponding Venn parts is just a permutation of the sets that contain them.See appendix 2 for illustration.

\item 1 unordered $k$-set tuple corresponds to 1 compatible unordered  $k$-event tuple:\\
The direction \emph{set}$\mapsto$ \emph{event} is trivial.The reverse direction \emph{event}$\mapsto$\emph{set} is guaranteed by the compatibility condition,since if 2 of the $k$ events correspond to the same set,then these two events are not distributionally compatible.
\end{itemize}
$\hfill\blacksquare$

Now we want to point out two crucial points concerning eq.(1).

First,it is an intuitive fact that for $k$ large enough,no $k$-event tuples will be distributionally compatible(i.e. A too large number of events can not be made simultaneously true.),and there exists a $k_{max}$ such that for $k>k_{max}$,the terms in eq.(1) are all zero(c.f.Theorem 2).It is also intuitive that we must have $k_{max}=k_{max}(n)\ll S_n\equiv\binom{n}{P_1}+\cdots+\binom{n}{P_t}$,so that this fact greatly reduces the computing power needed to find the exact values of the various Ramsey numbers.However,it is still extremely hard to find out the exact exact values of the Ramsey numbers from eq.(13) below.Maybe the order system introduced in appendix 1 can be of some help in this respect.

The precise value of $k_{max}$,the maximum number of events that are distributionally compatible,is hard to determine,although a rather crude upper bound is given in appendix 3.

Second,we know that for $P<r$,the binomial coefficient $\binom{P}{r}=0$.

By virtue of the $\binom{k}{2}$ compatibility conditions $P_{i_\mu
i_\nu}\le r-1$ ($\forall 1\le\mu<\nu\le k$)and the obvious fact which follows from Notation 3 that the $P$'s with more than 2 subscripts are no greater than the $P$'s with two 2 subscripts,we see that in the exponential $\mbox{txp}\left[\sum_{s=1}^{k}(-1)^s\mathop{\sum}_{1\leq\lambda_1<\cdots<\lambda_s\leq k}\binom{P_{i_{\lambda_1}\cdots i_{\lambda_s}}}{r}\right]$ in eq.(7),lots of terms in the sum $\mathop{\sum}_{1\leq\lambda_1<\cdots<\lambda_s\leq k}\binom{P_{i_{\lambda_1}\cdots i_{\lambda_s}}}{r}$ are actually zero. This simplification is automatically taken into account if we introduce the following variables.

Suppose that for a given unordered $k$-event tuple $(A_{i_1j_1},\cdots,A_{i_kj_k})$, $k_i$ $P_i$-events of the $k$ events demand their $r$-subsets be put into box $i$.($1\le k_i\le\binom{n}{P_i},\sum_{i=1}^{t}k_i=k$).Since the compatibility conditions amount to the statement that the events correspond to different boxes must have no common $r$-subsets.We can write
\begin{align}
&\left|\bigcup_{s=1}^{k}\mbox{Tr}(\bar{A}_{i_sj_s})\right|\notag\\
&=\sum_{i=1}^{t}\left|\bigcup_{\lambda=1}^{k_i}\mbox{Tr}(\bar{A}_{ij^{(i)}_\lambda})\right|\notag\\
&=\sum_{i=1}^{t}\sum_{s=1}^{k_i}(-)^s\sum_{1\le\lambda_1<\cdots<\lambda_s\le k_i}\left|\bigcap_{m=1}^{s}\mbox{Tr}(\bar{A}_{ij^{(i)}_{\lambda_m}})\right|\notag\\
&=\sum_{i=1}^{t}\sum_{s=1}^{k_i}(-)^s\sum_{1\le\lambda_1<\cdots<\lambda_s\le k_i}
\binom{\left|\bigcap_{m=1}^{s}\bar{A}_{ij^{(i)}_{\lambda_m}}\right|}{r}\notag\\
&=\sum_{i=1}^{t}\sum_{s=1}^{k_i}(-)^s\sum_{1\le\lambda_1<\cdots<\lambda_s\le k_i}
\binom{\sum\mbox{$2^{k-s}$\quad$Q$'s whose $(k_i+\lambda_m)^{th}$($1\le m\le s$) digits are 1}}{r}
\end{align}
where the unordered $k$-event tuple $(A_{i_1j_1}<\cdots<A_{i_kj_k})$ has been written as $A_{1j^{(1)}_1}<\cdots<A_{1j^{(1)}_{k_1}}<A_{2j^{(2)}_1}<\cdots<A_{2j^{(2)}_{k_2}}<\cdots<A_{tj^{(t)}_1}<\cdots<A_{tj^{(t)}_{k_t}}$.

The zero terms in the original sum $\mathop{\sum}_{1\leq\lambda_1<\cdots<\lambda_s\leq k}\binom{P_{i_{\lambda_1}\cdots i_{\lambda_s}}}{r}$ do not appear at all in the equation above .

Combining Theorems 2 and 3 and the two points mentioned above,we can summarize the main result of this paper as follows:
\begin{theorem}
The Ramsey number $R(P_1,P_2,...,P_t;r)$ is the smallest positive integer $n$ that satisfies
the following equation,($k\equiv k_1+\cdots k_t\le k_{max}(n)$)
\begin{align}
&\sum_{k_1=0}^{\binom{n}{P_1}}\cdots\sum_{k_t=0}^{\binom{n}{P_t}}(-)^{k-1}\sum_{\{Q\}} n!/\prod_{\mbox{all $2^k$\quad$B$'s}}Q^{(k)}_B!\times\notag\\
&\mbox{txp}\left[\sum_{i=1}^{t}\sum_{s=1}^{k_i}(-)^s\sum_{1\le\lambda_1<\cdots<\lambda_s\le k_i}
\binom{\sum\mbox{$2^{k-s}$\quad$Q$'s whose $(k_i+\lambda_m)^{th}$($1\le m\le s$) digits are 1}}{r}\right]\notag\\
&=1
\end{align}
,for some $k_{max}(n)\ll\binom{n}{P_1}+\cdots+\binom{n}{P_t}$.
Here,the summation $\sum_{\{Q\}}$ represents a $2^k$-fold conditional summation w.r.t the $2^k$ $Q^{(k)}$'s.The $Q^{(k)}$'s are constrained by the following $k+1$ conditions:
\begin{align}
&\sum(\mbox{all $2^k\quad Q^{(k)}$'s})=n \\
&\sum(\mbox{all $2^{k-1}\quad Q^{(k)}$'s whose $i_l^{th}$ digit is 1})=P_{i_l} (1\le l\le k)
\end{align}
together with the $\binom{k}{2}$ compatibility conditions of the unordered $k$-event tuple $(A_{i_1j_1},\cdots,A_{i_kj_k})$ stated in Def.2,translated from the language of the $P$'s to that of the $Q$'s via eq.(4):
\begin{align}
\sum(\mbox{all $2^{k-2}\quad Q^{(k)}$'s whose $i_\mu^{th}$ and $i_{\nu}^{th}$ digits are 1})\le r-1\quad(\forall 1\le\mu<\nu\le k)
\end{align}
,if the $Q$'s cannot be arranged to satisfy the compatibility conditions,then the term $\mbox{txp}\left[\cdots\right]$ should be replaced by zero.
\qed
\end{theorem}
\textbf{Proof}:
Theorems 2 and 3 offer the values of the two factors \emph{value of} $N(A_{i_1j_1}\bigcap\cdots\bigcap A_{i_kj_k})$ and \emph{frequency this value occurs} in eq.(6),respectively.What we need to consider here is the range of summation.The \emph{value} is a function of the Venn spectrum of the unordered $k$-event tuple $(A_{i_1 j_1},\cdots ,A_{i_k j_k})$.The essence of the summation $\sum_{\mbox{value}}$ and the second summation in eq.(1) is that we must traverse all the $\binom{S_n}{k}$ unordered $k$-event tuples $(A_{i_1 j_1},\cdots ,A_{i_k j_k})$,where $S_n\equiv\binom{n}{P_1}+\cdots+\binom{n}{P_t}$.The $t$-fold sum $\sum_{k_1=0}^{\binom{n}{P_1}}\cdots\sum_{k_t=0}^{\binom{n}{P_t}}$ exhausts all possible ``types'' $(k_1\cdots k_t)$ of the unordered $k$-event tuples that we are summing over.The third $\sum_{\{Q\}}$ traverses all distributionally compatible unordered $k$-event tuples $(A_{i_1 j_1},\cdots ,A_{i_k j_k})$ of this ``type''.The constraint (14) is a trivial constraint;the constraint (15) is the ``type'' constraint and the constraint (16)
is the ``compatibility'' constraint. $\hfill\blacksquare$
\section*{Appendix 1:Order systems}

We introduce order for the various structures encountered in this paper with the hope that unordered $k$-event tuples that are ``adjacent'' according to this order should have similar properties,so that this property(usually a function of unordered $k$-event tuples) can be calculated recursively from ``small'' unordered $k$-event tuples to ``big'' unordered $k$-event tuples.

Introducing order of unordered $k$-event tuples amounts to studying the geometry of a ``line'' in the space $V_{P_1,\cdots,P_t}^k$.(c.f. Defs.4 and 5)
We can associate other geometrical configurations with the space of unordered $k$-event tuples,w.r.t. which the relation among this property of different unordered $k$-event tuples may be more transparent.

The generic idea of order is as follows.

For two ``structures'' $X$ and $Y$ of the same type,we first compare their ``smallest'' ``elements'' $x_1$ and $y_1$.If $x_1=y_1$,then we compare $x_2$ and $y_2$...If finally we find some number $k$ such that $x_l=y_l,(\forall 1\le l<k)$ and $x_k<y_k$,then we say $X<Y$.

Using this idea ,all order structures ultimately reduce to the order structure of $V_n=\{1,\cdots,n\}$.

We begin to give definitions.All definitions are given in terms of events.The corresponding definitions in terms of sets are obvious.($\bigvee$ means OR,while $\bigwedge$ means AND.)
\begin{definition}
The $P_i$-event $\{a_1<\cdots<a_{P_i}\}$ is said to be ``less than'' the $P_i$-event $\{b_1<\cdots<b_{P_i}\}$:\\
$\{a_1<\cdots<a_{P_i}\}<\{b_1<\cdots<b_{P_i}\}$
iff $\exists s,\forall l<s, a_l=b_l \bigwedge a_s<b_s$.

When denoting them by $A_{ij_1}=\{a_1<\cdots<a_{P_i}\}$ and $A_{ij_2}=\{b_1<\cdots<b_{P_i}\}$,we should assign the labels
$j_1<j_2$.
\qed
\end{definition}

\begin{definition}
The $P_i$-event $A_{ij}$ is said to be ``less than'' the $P_l$-event $A_{lm}$:\\
$A_{ij}<A_{lm}$
iff $(i<l)\bigvee[(i=l)\bigwedge(j<m)]$
\qed
\end{definition}
\begin{definition}
The unordered $k$-event tuple
$(A_{i_1j_1},\cdots,A_{i_kj_k})$ is said to be ``less than'' the unordered $m$-event tuple $(A_{u_1v_1},\cdots,A_{u_mv_m})$:\\
$(A_{i_1j_1},\cdots,A_{i_kj_k})<(A_{u_1v_1},\cdots,A_{u_mv_m})$
 iff
$(k<m)\bigvee\{(k=m)\bigwedge\exists s[\forall l<s ,A_{i_lj_l}=A_{u_lv_l}]\bigwedge[A_{i_sj_s}<A_{u_sv_s}]\}$
\qed
\end{definition}
\section*{Appendix 2:Three-set Venn diagram that illustrates one part in the proof of Theorem 3}

\begin{figure}[!ht]
        \centering
                \includegraphics[scale=1.7]{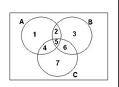}
                       \caption{3-set Venn diagram}
        \label{Fig:1}
\end{figure}

In Fig.1,the $n$-set $V_n$ is divided into 8 Venn parts labeled by 1,2,...,7.(The Venn part $A^c\bigcap B^c\bigcap C^c$ is not labeled by the number 8.)

We can \emph{simultaneously} interchange the parts $1\leftrightarrow3$ and $4\leftrightarrow6$.This generates a \emph{new} way of partition of the $n$-set $V_n$ into 8 \emph{ordered} Venn parts,but correspond to the \emph{same} \emph{unordered} 3-set tuple $(A,B,C)$,since the net effect of this interchange is just an exchange of $A$ and $B$, which is one of the 3! permutations of the 3 sets $A,B$ and $C$ that do not change the unordered 3-event tuple $(A,B,C)$.

There are also other simultaneous permutations of the corresponding Venn parts that do not change the unordered 3-event tuple $(A,B,C)$.For instance, the two simultaneous permutations $1\rightarrow7\rightarrow3\rightarrow1$ and $2\rightarrow4\rightarrow6\rightarrow2$ have the effect $A\rightarrow C\rightarrow B\rightarrow A;(A,B,C)=(B,C,A)$,while $1\rightarrow3\rightarrow7\rightarrow1$ and $2\rightarrow6\rightarrow4\rightarrow2$ have the effect $A\rightarrow B\rightarrow C\rightarrow A;(A,B,C)=(C,A,B)$.
\section*{Appendix 3:An upper bound for $k_{max}$}
The number of $P_j$-events that are distributionally compatible with a given $P_i$-event($j\ne i$) is
$\sum_{\nu=0}^{r-1}\binom{P_i}{\nu}\binom{n-P_i}{P_j-\nu}$,which is just the sum of the number of $P_j$-events events that have $\nu$ common elements with the given $P_i$-event.($\le\nu\le r-1$).

Therefore if this $P_i$-event is chosen as one of the $k$ events in a given $k$-event tuple, then the remaining $k-1$ events can only be chosen from the $\sum_{j=1,j\ne i}^{t}\sum_{\nu=0}^{r-1}\binom{P_i}{\nu}\binom{n-P_i}{P_j-\nu}+\binom{n}{P_i}-1$ events that are compatible with this $P_i$-event.

Since the above argument is valid for any $i$,we have for $k_{max}$, the maximum number of events that are distributionally compatible,the upper bound
\begin{align}
k_{max}(n)\le\max_i\left\{\binom{n}{P_i}+\sum_{j=1,j\ne i}^{t}\sum_{\nu=0}^{r-1}\binom{P_i}{\nu}\binom{n-P_i}{P_j-\nu}\right\}
\end{align}


\begin{thebibliography}{99}
\bibitem{ram} J.H.Van Lint and R.M.Wilson \emph{A course in Combinatorics} second edition Cambridge
University Press (2001)
\end{thebibliography}
\end{document}